\documentclass[a4paper,12pt]{article}
\usepackage{amsmath}
\usepackage{amssymb}

\begin{document}
\title{\bf On the Sendov conjecture for polynomials with simple zeros}
\author{{\bf Robert DALMASSO}\\\\ Le Galion - B,\\ 
33 Boulevard Stalingrad,\\
 06300 Nice, France}

\date{}
\maketitle

\allowdisplaybreaks[4]

\bigskip

\hrule

\bigskip

\noindent{\bf Abstract.} The Sendov conjecture asserts that if all the zeros of a polynomial $p$ lie in the closed unit disk then there must be a zero of $p'$ within unit distance of each zero. In this paper we give a partial result when $p$ has simple zeros.

\smallskip


\smallskip


\bigskip

\hrule

\bigskip

\section{Introduction}

The well-known Sendov conjecture (\cite{h} Problem 4.5) asserts that if $p(z) = \prod_{j=1}^{n}(z - z_j)$ is a polynomial with $|z_j| \leq 1$ ($1 \leq j \leq n$), then each disk $|z - z_j|\leq 1$ ($1 \leq j \leq n$) contains a zero of $p'$. Notice that by the Gauss-Lucas theorem the zeros $w_k$ ($1 \leq k \leq n-1$) of $p'$ lie in the closed convex hull of the zeros of $p$, hence $|w_k| \leq 1$ for $1 \leq k \leq n-1$. 

This conjecture has been verified for polynomials of degree $n \leq 8$ or for arbitrary degree $n$ if there are at most eight distinct roots: See Brown and Xiang \cite{bx} and the references therein. It is also true in general ($n \geq 2$)  when $p(0) = 0$ (\cite{sc}). The Sendov conjecture is true with respect to the root $z_j$ of $p$ if $|z_j| = 1$ (\cite{ru}). Recently (\cite{de}) it has been verified when $n$ is larger than a fixed integer depending on the root $z_j$ of $p$. We refer the reader to  Marden \cite{m} and  Sendov \cite{s} for further information and bibliographies.

\medskip

In this paper we prove the following theorem.

\bigskip

\noindent{\bf Theorem 1} {\sl Suppose that $p$ has simple zeros $z_j$, $j = 1,\cdots, n$ in the closed unit disk $\overline D(0,1)$. Let
\begin{displaymath}
r = \frac{1}{n} \min_{1\leq j \leq n}\min_{i\not= j}|z_i - z_j| \quad , \, \,  \displaystyle a_n = n^{\frac{1}{n-1}} - 1 \, ,
\end{displaymath}
\begin{displaymath}
A_{2p} = 2p(\frac{r^{2}a_{2p}}{24})^p  \quad \textrm{and} \, \,  \displaystyle \, \, A_{2p+1} = (2p+1)a_{2p+1}^{p+1}(\frac{r^{2}\sqrt3}{72})^p\, .
\end{displaymath}
Then the Sendov conjecture holds for $p$ when $|p(0)| \leq A_n$.}

\bigskip

Our proof will make use of the the so-called Coincidence Theorem, a variant of Grace's Apolarity Theorem (\cite{m}). We begin with the following definition.

\bigskip

\noindent{\bf Definition} {\sl $\Phi(x_1,\cdots,x_p)$ is a symmetric $p$-linear form of total degree $p$ in the variables $x_j$ ( $1\leq j \leq p$) belonging to $\mathbb C$ if it is symmetric in these variables and if $\Phi$ is a polynomial of degree 1 in each $x_j$ separately such that $\Phi(x,\cdots,x)$ is a polynomial of degree $p$ in $x$.}

\bigskip

A {\sl circular region} is an open or closed disk or halfplane in $\mathbb C$, or the complement of any such set. Now we recall Walsh's Coincidence Theorem (\cite{m}).

\bigskip

\noindent{\bf Theorem 2} {\sl Let $\Phi$ be a symmetric  $p$-linear form of total degree $p$ in $x_1,\cdots,x_p$ and let $C$ be a circular region containing the $p$ points $\alpha_1,\cdots,\alpha_p$. Then in $C$ there exists at least one point $x$ such that $\Phi(\alpha_1,\cdots, \alpha_p) = \Phi(x,\cdots,x)$.}

\bigskip

Finally we also need the next result.

\bigskip

\noindent{\bf Theorem 3} (\cite{di}) {\sl Let 
\begin{displaymath}
q(z) = \prod_{j=1}^{r}(z-u_j)^{k_j} \, \, , \quad \sum_{j=1}^{r} k_j = m,
\end{displaymath}
be a polynomial of degree $m$ whose zeros $u_1 ,\cdots , u_r$ are distinct and have multiplicity $k_1,\cdots, k_r$ , respectively. For any zero $u_j$ of $q$, let $M_j = \min_{i\not= j}|u_i-u_j|$, $j = 1,\cdots, r$. Then $q$ has no nontrivial critical point (the critical points which are not zeros of the polynomial) in
\begin{displaymath}
\displaystyle \bigcup_{j=1}^{r}\{z \in \mathbb C \, ; |z-u_j| < \frac{k_j}{m}M_j\, \} \, .
\end{displaymath}}

\medskip

In Section 2 we give some preliminary results. Theorem 1 is proved in Section 3.

\section{Preliminaries}

We begin with the following lemma.

\bigskip

\noindent{\bf Lemma 1} {\sl The Sendov conjecture is true with respect to the root $z_j$ if $|z_j| \leq a_n$.}

\bigskip

\noindent{\sl Proof.} Since $p'(z_j) = q(z_j)$, where $q(z) = p(z)/(z-z_j)$, we have
\begin{displaymath}
n\prod_{k=1}^{n-1}(z_j - w_k) = \prod_{k \not= j}(z_j - z_k) \,.
\end{displaymath}
Then
\begin{displaymath}
n\prod_{k=1}^{n-1}|z_j - w_k| = \prod_{k \not= j}|z_j - z_k| \leq (1+|z_j|)^{n-1} \leq n \,,
\end{displaymath}
and the lemma follows.

\medskip

Distinguish one of the zeros of $p$, say $z_n$, and let $z_n = a$. By a rotation, if necessary, and using Lemma 1 we may suppose that $a_n < a < 1$. Let $0 < s \leq r$ and set $z_0 = \displaystyle\frac{a}{2} + \frac{i}{2}\sqrt{4-a^2}$. Let $v_1(z_0,s)$, $v_2(z_0,s) \in \partial D(0,1)\cap\partial D(z_0,s)$  with $ \textrm{Re}\, v_1(z_0,s) <  \textrm{Re}\, v_2(z_0,s)$. We denote by $L(z_0,s)$ the line through $v_1(z_0,s)$ and $v_2(z_0,s)$ and by $H(z_0,s)$ the closed halfplane bounded by $L(z_0,s)$ such that $0 \notin H(z_0,s)$. We set $A(z_0,s) = \overline D(0,1) \cap H(z_0,s)$. Finally let $v_3(z_0,s) \in L(z_0,s)\cap\partial D(a,1)$ with $ \textrm{Re}\, v_3(z_0,s) < a/2$.

\bigskip

\noindent{\bf Lemma 2} {\sl With the above notations we have
\begin{displaymath}
{\rm Re}\, v_1(z_0,s) = \frac{1}{4}(a(2-s^2) - s((4-a^2)(4-s^2))^{\frac{1}{2}}) \,,
\end{displaymath}
\begin{displaymath}
{\rm Im}\, v_1(z_0,s) = \frac{-a}{4\sqrt{4-a^2}}(a(2-s^2) - s((4-a^2)(4-s^2))^{\frac{1}{2}}) + \frac{2-s^2}{\sqrt{4-a^2}} \,,
\end{displaymath}
\begin{displaymath}
{\rm Re}\, v_2(z_0,s) = \frac{1}{4}(a(2-s^2) + s((4-a^2)(4-s^2))^{\frac{1}{2}}) \,,
\end{displaymath}
\begin{displaymath}
{\rm Im}\, v_2(z_0,s) = \frac{-a}{4\sqrt{4-a^2}}(a(2-s^2) + s((4-a^2)(4-s^2))^{\frac{1}{2}}) + \frac{2-s^2}{\sqrt{4-a^2}} \,,
\end{displaymath}
\begin{displaymath}
{\rm L(z_0,s)}\, = \{ c + id \, ; \, c\, , \, d \in \mathbb R \quad {\rm and} \quad d\sqrt{4-a^2} = - ac + 2-s^2 \} \,,
\end{displaymath}
and
\begin{displaymath}
{\rm Re}\, v_3(z_0,s) = \frac{1}{4}(a(6-a^2 -s^2) - ((4-a^2)(4-a^2 -s^2)(a^2 + s^2))^{\frac{1}{2}}) \, .
\end{displaymath}
Moreover 
\begin{displaymath}
\frac{a}{2} - {\rm Re}\, v_3(z_0,s) > \frac{s^2}{6} \, .
\end{displaymath}
}

\bigskip

\noindent{\sl Proof.} We only prove the  inequality since all the other formulas follow from elementary computations. Set
\begin{displaymath}
\Delta = (4-a^2)(4-a^2 -s^2)(a^2+s^2) \, .
\end{displaymath}
We can write
\begin{displaymath}
\begin{array}{lcl}
\displaystyle\frac{a}{2} - {\rm Re}\, v_3(z_0,s)
& = &\displaystyle \frac{a}{2} + \frac{1}{4}(\sqrt{\Delta} - a(6-a^2-s^2))  \\ \\
& = &\displaystyle \frac{a}{2} + \frac{1}{4}\frac{\Delta - a^2(6-a^2-s^2)^2}{\sqrt{\Delta} + a(6-a^2-s^2)}    \\ \\
& = &  \displaystyle \frac{1}{4}\frac{2a\sqrt{\Delta} + 2a^2(6-a^2-s^2) + \Delta - a^2(6-a^2-s^2)^2}{\sqrt{\Delta} + a(6-a^2-s^2)} \\ \\ 
& = & \displaystyle\frac{1}{4}\frac{N}{D}\, .
\end{array}
\end{displaymath}
Since
\begin{displaymath}
\sqrt{\Delta} > \displaystyle a((4-a^2)(4-a^2 -s^2))^{\frac{1}{2}} > a(4-a^2-s^2) \, ,
\end{displaymath}
we get
\begin{displaymath}
\begin{array}{lcl}
N
& = &\displaystyle 2a\sqrt{\Delta} + 2a^2(6-a^2-s^2) + \Delta - a^2(6-a^2-s^2)^2  \\ \\
& = &\displaystyle 2a\sqrt{\Delta} + \Delta   - a^2(6-a^2-s^2)(4-a^2-s^2) \\ \\ 
& > & \displaystyle 2a^2(4-a^2-s^2) + s^2(4-a^2)(4-a^2 -s^2) \\ \\
& & - a^2(4-a^2-s^2)(2-s^2) \\ \\
& = & 4s^2(4-a^2-s^2) > 8s^2 \, .
\end{array}
\end{displaymath}
Since $D < \displaystyle\sqrt{32} + 6 < 12$, we get
\begin{displaymath}
\frac{a}{2} - {\rm Re}\, v_3(z_0,s) > \frac{s^2}{6} \, .
\end{displaymath}

\noindent{\bf Lemma 3} {\sl Suppose that there exists $k \in \{1, \cdots, n-1\}$ such  that ${\rm Re}\, w_k \leq a/2$ and ${\rm Im}\, w_k \geq {\rm Im}\, \displaystyle v_3(z_0,\frac{r}{2})$. Then ${\rm Re}\, w_k < \displaystyle{\rm Re}\,\displaystyle v_3(z_0,\frac{r}{2})$.
}

\bigskip

\noindent{\sl Proof.} We claim that $p'$ has no critical point in $\displaystyle A(z_0,\frac{r}{2})$. Then the lemma follows. To verify this claim suppose first that $p(z_0) = 0$. Since $p$ has simple zeros Theorem 3 implies that $p'(z) \not= 0$ for all $z \in D(z_0 ,r)$ and the claim is true in this case. Now assume that $p(z_0) \not= 0$. If there exists $z_j \in \displaystyle \overline D(z_0,\frac{r}{2})\cap\overline D(0,1)$ such that $p(z_j) = 0$  for some $j \in \{1,\cdots, n- 1\}$, again Theorem 3 implies that $p'$ has no critical point in $D(z_j,r)$ and the claim is true. Finally, if  $p(z) \not= 0$  for  $z \in \displaystyle  \overline D(z_0,\frac{r}{2})\cap\overline D(0,1)$, the Gauss-Lucas Theorem implies that $p'(z) \not= 0$ for $\displaystyle z \in A(z_0,\frac{r}{2})$. 
\bigskip

\noindent{\sl Remark 1.} Let $z'_0 = \displaystyle \frac{a}{2} - \frac{i}{2}\sqrt{4-a^2}$. Lemmas similar to Lemma 2 and Lemma 3 hold if we consider the point $z'_0$ instead of $z_0$.

\bigskip

\noindent{\bf Lemma 4} {\sl Let
\begin{displaymath}
\alpha_j = \frac{\sin(2j\pi)/n}{1 - \cos(2j\pi)/n} \quad , \,\, 1 \leq j \leq n-1 \, .
\end{displaymath}
Then $\alpha_1\cdots\alpha_{p-1} = 1$ when $n = 2p$ and  $\alpha_1\cdots\alpha_p \geq (\sqrt3/3)^{p}$ when $n = 2p+1$.}

\bigskip

\noindent{\sl Proof.} Let $n = 2p$. For $j \in \{1,\cdots,p-1\}$ we have
\begin{displaymath}
\displaystyle \alpha_{p-j} =   \frac{\sin(j\pi/p)}{1+\cos(j\pi/p)} = \frac{1 - \cos(j\pi/p)}{\sin(j\pi/p)} =\frac{1}{\alpha_j}\quad ,
\end{displaymath}
and the result follows. 

Now let $n = 2p+1$ with $p \geq 2$. We can write
\begin{displaymath}
\alpha_1\cdots\alpha_p  = \prod_{j=1}^{k}\alpha_j\alpha_{p-j+1} \quad \textrm{if} \, \, p = 2k, \, \, k \geq 1 \,, 
\end{displaymath}
and
\begin{displaymath}
\alpha_1\cdots\alpha_p  = \alpha_{k+1}\prod_{j=1}^{k}\alpha_j\alpha_{p-j+1} \quad \textrm{if} \, \, p = 2k+1, \, \, k \geq 1 \, . 
\end{displaymath}
For $j \in \{1,\cdots,k\}$ we have
\begin{displaymath}
\begin{array}{lcl}
\displaystyle \alpha_j\alpha_{p-j+1}
& = & \displaystyle \frac{\sin(2j\pi/(2p+1))}{1 - \cos(2j\pi/(2p+1))} \frac{\sin((2j-1)\pi/(2p+1))}{1 + \cos((2j-1)\pi/(2p+1))}\\ \\
& = & \displaystyle \frac{\cos(j\pi/(2p+1))}{\sin(j\pi/(2p+1))} \frac{\sin((2j-1)\pi/(2p+1))}{1 + \cos((2j-1)\pi/(2p+1))}\\ \\
& \geq & \displaystyle \frac{\cos(j\pi/(2p+1))}{\sin(j\pi/(2p+1))} \frac{\sin(j\pi/(2p+1))}{1 + \cos(j\pi/(2p+1))}\\ \\
& = &\displaystyle \frac{\cos(j\pi/(2p+1))}{1 + \cos(j\pi/(2p+1))}  \geq \frac{\cos(\pi/3)}{1 + \cos(\pi/3)}  = \frac{1}{3} \, .
\end{array}
\end{displaymath}
Now when $p = 2k+1$ we have
\begin{displaymath}
\displaystyle \alpha_{k+1} = \frac{\cos((k+1)\pi/(4k+3))}{\sin((k+1)\pi/(4k+3))} \geq \frac{\cos(\pi/3)}{\sin(\pi/3)}  = \frac{\sqrt3}{3} \, .
\end{displaymath}

The  lemma follows.

\bigskip

\noindent{\bf Lemma 5} {\sl 1) We have

\begin{displaymath}
\frac{r^2}{24} < \frac{a_{2p}}{2}\alpha_{p-1} \quad \textrm{if} \, \, n =2p \quad \textrm{and} \quad 
\frac{r^2}{24} < \frac{a_{2p+1}}{2}\alpha_{p} \quad \textrm{if} \, \, n =2p +1 \, .
\end{displaymath}

2) Let $v \in \mathbb C$ be such that
\begin{displaymath}
{\rm Re}\,  v  < \frac{a}{2} -  b \, ,
\end{displaymath}
where $b > 0$ is such that
\begin{displaymath}
b < \frac{a_{2p}}{2}\alpha_{p-1} \quad \textrm{if} \, \, n =2p \quad \textrm{and} \quad 
b < \frac{a_{2p+1}}{2}\alpha_{p} \quad \textrm{if} \, \, n =2p +1 \, .
\end{displaymath}
Then 
\begin{displaymath}
|(v-a)^n - v^n| >  B_n(b) \, ,
\end{displaymath}
where
\begin{displaymath}
B_{2p}(b) = 2p(ba_{2p})^p  \quad \textrm{and} \, \,  \displaystyle \, \, B_{2p+1} = (2p+1)a_{2p+1}^{p+1}(b\frac{\sqrt3}{3})^p\, .
\end{displaymath}}

\bigskip

\noindent{\sl Proof.} 1) is easily verified since $r \leq 2/n$.

2) Let 
\begin{displaymath}
C_n(y) = \prod_{j=1}^{n-1}(b^2 + (y-\frac{a}{2}\alpha_j)^2) \quad, \,\, y \in \mathbb R \,\, .
\end{displaymath}
We have
\begin{displaymath}
C_{2p}(y) = (b^2 + y^2)\prod_{j=1}^{p-1}(b^2 + (y-\frac{a}{2}\alpha_j)^2)(b^2 + (y+\frac{a}{2}\alpha_j)^2)\quad , 
\end{displaymath}
and
\begin{displaymath}
C_{2p+1}(y) = \prod_{j=1}^{p}(b^2 + (y-\frac{a}{2}\alpha_j)^2)(b^2 + (y+\frac{a}{2}\alpha_j)^2) \,\, .
\end{displaymath}
Since $a > a_n$ we get
\begin{displaymath}
(b^2 + (y-\frac{a}{2}\alpha_j)^2)(b^2 + (y+\frac{a}{2}\alpha_j)^2) \geq b^2a^2\alpha_j^2 \quad  \textrm{for} \, \, 1 \leq j \leq p \quad  \textrm{and} \, \, y \in \mathbb R\, .
\end{displaymath}
Then Lemma 4 implies that
\begin{displaymath}
C_{2p}(y) >  b^{2p}a_{2p}^{2(p-1)} \quad  \textrm{and} \quad C_{2p+1}(y) > b^{2p}a_{2p+1}^{2p}3^{-p}  \quad  \textrm{for} \, \, y \in \mathbb R\, .
\end{displaymath}
Now the solutions $v_j$ ($1 \leq j \leq n-1$) of 
\begin{displaymath}
(v - a)^n - v^n = 0 \,\, ,
\end{displaymath}
are given by\begin{equation}
 v_j = \frac{a}{2}(1 + i \alpha_j) \, \, , \, \, 1 \leq j \leq n-1 \, .
 \label{eq:eq1}
 \end{equation}
Therefore 
\begin{displaymath}
\begin{array}{lcl}
\displaystyle |(v-a)^n - v^n|^2 & = & \displaystyle n^2a^2\prod_{j=1}^{n-1}(({\rm Re}\, v -\frac{a}{2})^2 + ({\rm Im}\, v - \frac{a}{2}\alpha_j)^2) \\ \\
& \geq & n^2a^2C_n({\rm Im}\, v) \\ \\
& > &  B_n(b)^2\, .
\end{array}
\end{displaymath}

\section{Proof of Theorem 1}

Define
\begin{displaymath}
\begin{array}{lcl}
k(z,u_1,\cdots,u_{n-1}) & = & z^n + \displaystyle \sum_{k=1}^{n-1}(-1)^{k}\frac{n}{n-k}\\ \\
&  & \displaystyle \times(\sum_{1\leq i_1 < \cdots < i_k \leq n-1}u_{i_1}\cdots u_{i_k})z^{n-k} \, \, ,
\end{array}
\end{displaymath}  
for $z, u_1, \cdots,u_{n-1} \in \mathbb C$.
We have 
\begin{displaymath}
k(z,w_1,\cdots,w_{n-1}) = p(z) - p(0) \, \, .
\end{displaymath}  
Let
\begin{displaymath}
C = \{ z \in \mathbb C \, ; \, {\rm Re}\, z <  \displaystyle \frac{a}{2} - \frac{r^2}{24} \, \} \, .
\end{displaymath}
Suppose that $|a-w_k| > 1$ for $k = 1,\cdots,n-1$. Then ${\rm Re}\, w_k < a/2$ for $k = 1,\cdots,n-1$. Using Lemma 2 with $s = r/2$, Lemma 3 and Remark 1  we obtain ${\rm Re}\, w_k < \displaystyle \frac{a}{2} - \frac{r^2}{24}$ for $k = 1,\cdots,n-1$. Theorem 2 implies that there exists $v \in C$ such that 
\begin{displaymath}
- p(0) = k(a,w_1,\cdots,w_{n-1}) = k(a,v,\cdots, v) = (a - v)^n + (-1)^{n-1}v^n \, .
\end{displaymath}
Lemma 5  with $b = r^2/24$ implies that  $v \notin C$ and we reach a contradiction.

\bigskip

\noindent{\sl Remark 2.} Suppose that the zeros of $p$ are not necessarily simple. Using Walsh's Coincidence Theorem we can give a new proof of Sendov's conjecture when $p(0) = 0$.
By a rotation, if necessary, we may suppose that $p$ has the form
\begin{displaymath}
p(z) = (z-a)\prod_{j=1}^{n-1}(z-z_j) \, ,
\end{displaymath}
where $a \in (0,1]$ is a simple root. Suppose that $|a-w_k| > 1$ for $k = 1, \cdots , n-1$. Since $ \textrm{Re}\, w_k < a/2$, there exists $t > 0$ such that $ \textrm{Re}\, w_k < \displaystyle \frac{a}{2} - t$ for $k = 1, \cdots , n-1$. Let $E = \{z \in \mathbb C \, ;\, \textrm{Re}\, z < \displaystyle \frac{a}{2} - t\,\}$. Theorem 2 implies that there exists $v \in E$ such that 
\begin{displaymath}
0 = k(a,w_1,\cdots,w_{n-1}) = k(a,v,\cdots, v) = (a - v)^n + (-1)^{n-1}v^n \, .
\end{displaymath}
Since $v = v_j$ for some $j \in \{1,\cdots, n-1\}$, \eqref{eq:eq1} implies that $v \notin E$ and we reach a contradiction.

\end{document}